\documentclass[12pt]{amsart}
\usepackage{pgfplots}
\usepackage{tikz}

\usetikzlibrary{automata}
\usetikzlibrary{positioning}
\usepackage{amsmath, amssymb} 
\usepackage{amsthm}
\usepackage{graphicx}
\usepackage[shortlabels]{enumitem}

\usepackage{amsfonts,geometry,hyperref}

\newtheorem{conj}{Conjecture}[section]

\newtheorem{theorem}{Theorem}[section]

\newtheorem{lemma}[theorem]{Lemma}

\newtheorem{prob}[theorem]{Problem}

\DeclareMathOperator{\Nbd}{Nbd}
\DeclareMathOperator{\ABC}{ABC}
\DeclareMathOperator{\GA}{GA}
\DeclareMathOperator{\SCI}{SCI}

\begin{document}
\title[Zagreb indices of subgroup generating bipartite graph]{Zagreb indices of subgroup generating bipartite graph}

\author[S. Das, A. Erfanian and R. K. Nath]{Shrabani Das, Ahmad Erfanian and Rajat Kanti Nath$^*$}

\address{S. Das, Department of Mathematical Sciences, Tezpur University, Napaam-784028, Sonitpur, Assam, India.}

\email{shrabanidas904@gmail.com}

\address{A. Erfanian, Department of Pure Mathematics, Ferdowsi University of Mashhad,
	P.O. Box 1159-91775,
	Mashhad, Iran}
	\email{erfanian@um.ac.ir}

\address{R. K. Nath, Department of Mathematical Sciences, Tezpur University, Napaam-784028, Sonitpur, Assam, India.} 
\email{ rajatkantinath@yahoo.com}
\thanks{$^*$Corresponding Author}
\begin{abstract}
Let $G$ be a group and $L(G)$ be the set of all subgroups of $G$. The subgroup generating bipartite graph $\mathcal{B}(G)$ defined on $G$ is a bipartite graph whose vertex set is partitioned into two sets $G \times G$ and $L(G)$, and two vertices $(a, b) \in G \times G$ and $H \in L(G)$ are adjacent if $H$ is generated by $a$ and $b$. In this paper, we deduce expressions for first and second Zagreb indices of $\mathcal{B}(G)$ and  obtain a condition such that $\mathcal{B}(G)$ satisfy Hansen-Vuki{\v{c}}evi{\'c} conjecture [Hansen, P. and Vuki{\v{c}}evi{\'c}, D. Comparing the Zagreb indices, {\em Croatica Chemica Acta}, \textbf{80}(2), 165-168, 2007]. It is  shown that $\mathcal{B}(G)$  satisfies   Hansen-Vuki{\v{c}}evi{\'c} conjecture if $G$ is a cyclic group of order $2p, 2p^2, 4p$, $4p^2$ and $p^n$; dihedral group of order $2p$ and $2p^2$; and dicyclic group of order $4p$ and $4p^2$ for any prime $p$.   
While computing Zagreb indices of $\mathcal{B}(G)$ we have computed $\deg_{\mathcal{B}(G)}(H)$ for all $H \in L(G)$ for the above mentioned groups. Using these information we also compute Randic Connectivity index, Atom-Bond Connectivity index, Geometric-Arithmetic index, Harmonic index and Sum-Connectivity index  of $\mathcal{B}(G)$.
\end{abstract}

\thanks{ }

\maketitle

\section{Introduction}
Defining graphs on finite groups and studying them towards the characterization of finite groups/graphs is an old practice. Cayley graph, generating graph, power graph, commuting graph, nilpotent graph, solvable graph, commuting/ nilpotent/ solvable conjugacy class graph, B superA graph etc. (see \cite{cameron2021graphs,BNN2020,CaJa-2024,ANC-2022}) are  some well-known examples of such graphs. Note that all these graphs are not bipartite.
Let $G$ be a finite group and $L(G) := \{H : H \text{ is a subgroup of } G\}$. In \cite{DEN-23}, we had introduced a bipartite graph $\mathcal{B}(G)$  whose vertex set $V(\mathcal{B}(G))$ is partitioned into two sets $G \times G$ and $L(G)$, and two vertices $(a, b) \in G \times G$ and $H \in L(G)$ are adjacent if  $H = \langle a, b \rangle$, the subgroup generated by $a$ and $b$. We call this graph the subgroup generating bipartite graph (abbreviated as SGB-graph). We had obtained relations between $\mathcal{B}(G)$  and various probabilities (such as probability generating  a given subgroup \cite{Di69}, commuting probability \cite{DNP-13}, cyclicity degree \cite{PSSW93}, nilpotency degree \cite{DGMW92}, solvability degree \cite{FGSV2000}) associated to finite groups. The probability generating  a given subgroup $H$ of $G$, denoted by $\Pr_H(G)$, is the probability that a randomly chosen pair of elements of $G$ generate $H$. The origin of $\Pr_H(G)$ lies in a paper of Hall \cite{Hall36} and a generalized version of which was studied in \cite{Pak99}.
For any subgroup $H$ of $G$, it was shown  (see \cite[Lemma 3.1]{DEN-23}) that
\begin{equation}\label{deg(H in L(G))}
	\deg_{\mathcal{B}(G)}(H)=|G|^2 {\Pr}_H(G),
\end{equation}
where $\deg_{\mathcal{G}}(x)$ is the degree of any vertex $x$ in any graph $\mathcal{G}$. In \cite{DEN-23}, we had also 
discussed about various graph parameters of $\mathcal{B}(G)$ such as independence number, domination number, girth, diameter etc. 
It is worth mentioning that Mahtabi et al. \cite{MEM-2023} also considered a bipartite graph on $G$ using the automorphism group of $G$.  

 Let $\Gamma$ be the set of all graphs. A topological index  is a function $T :  \Gamma \to \mathbb{R}$ such that $T(\mathcal{G}_1) = T(\mathcal{G}_2)$ whenever the graphs $\mathcal{G}_1$ and   $\mathcal{G}_2$ are isomorphic. Topological indices of graphs are numerical invariants which provide information about their physical and chemical properties in case of chemical graphs. Many topological indices have been defined since 1947 by using different parameters of graphs. Some of them have proved to be useful also in other areas where connectivity patterns play an important role, such as interprocessor connections and complex networks.  One such popular topological indices are the Zagreb indices which are degree-based topological indices. They were introduced by Gutman and Trinajsti{\'c} \cite{Gut-Trin-72} in 1972. Zagreb index is used in examining the dependence of total $\pi$-electron energy on molecular structure. As noted in \cite{Z-index-30y-2003}, Zagreb index is also used in studying molecular complexity, chirality, ZE-isomerism and heterosystems etc. A survey on mathematical properties of Zagreb index can be found in \cite{Gut-Das-2004}. Zagreb index for commuting and non-commuting graphs of finite non-abelian groups have been studied in \cite{DSN-23, mirzargar2012some}. Zagreb indices of super commuting graphs of finite groups  have been computed in \cite{DN-24}.

  Let $\mathcal{G}$ be a simple undirected graph with vertex set $V(\mathcal{G})$ and edge set $e(\mathcal{G})$. The first and second Zagreb indices of $\mathcal{G}$, denoted by $M_{1}(\mathcal{G})$ and $M_{2}(\mathcal{G})$ respectively, are defined as 
\[
M_{1}(\mathcal{G}) = \sum\limits_{v \in V(\mathcal{G})} \deg(v)^{2}  \text{ and }  M_{2}(\mathcal{G}) = \sum\limits_{uv \in e(\mathcal{G})} \deg(u)\deg(v).
\]
Comparing first and second Zagreb indices, Hansen and Vuki{\v{c}}evi{\'c} \cite{hansen2007comparing} posed the following   conjecture in  2007.
\begin{conj}\label{Conj}
	(Hansen-Vuki{\v{c}}evi{\'c} Conjecture) For any simple finite graph $\mathcal{G}$, 
	\begin{equation}\label{Conj-eq}
		\dfrac{M_{2}(\mathcal{G})}{\vert e(\mathcal{G}) \vert} \geq \dfrac{M_{1}(\mathcal{G})}{\vert V(\mathcal{G}) \vert} .
	\end{equation}
\end{conj}
\noindent They disproved the conjecture in the same paper by showing that it is not true  if $\Gamma = K_{1, 5} \sqcup K_3$, where $\sqcup$ stands for disjoint union of graphs or sets. However, Hansen and Vuki{\v{c}}evi{\'c} \cite{hansen2007comparing}  showed that  Conjecture \ref{Conj}  holds for chemical graphs. 
Upon further studies, it was shown that the conjecture holds for trees and with equality in \eqref{Conj-eq} when $\Gamma$ is a star graph \cite{vukicevic2007comparing}. Thereafter, it was shown that the conjecture holds for connected unicyclic graphs with equality when the graph is a cycle (see \cite{liu2008conjecture}). The case when equality holds in \eqref{Conj-eq} is studied extensively in \cite{vukicevic2011some}. A survey on comparing Zagreb indices can be found in \cite{Liu-You-2011}. Zagreb indices of commuting conjugacy class graph are computed and verified  Conjecture \ref{Conj} in \cite{Das-Nath-2023} for the classes of finite groups considered in \cite{Salah-2020, SA-2020, SA-CA-2020}. 
 
In Section 2, we discuss about the structures of  $\mathcal{B}(G)$ for certain groups. In Section 3, we obtain expressions for Zagreb indices of $\mathcal{B}(G)$ in relation with $\Pr_H(G)$ and obtain a condition such that $\mathcal{B}(G)$ satisfy Hansen-Vuki{\v{c}}evi{\'c} conjecture. We also show that $\mathcal{B}(D_{2n})$ and $\mathcal{B}(Q_{4n})$  satisfies Hansen-Vuki{\v{c}}evi{\'c} conjecture if $n = p$ and $p^2$, where $D_{2n}$ is the dihedral group presented by $\langle a, b: a^n=b^2=1, bab=a^{-1} \rangle$, $Q_{4n}$ is the dicyclic group presented by $Q_{4n} =  \langle a, b : a^{2n} = 1, b^2 = a^n, bab^{-1} = a^{-1} \rangle$ and $p$ is a prime. For a cyclic group $G$ and a prime $p$, we show that $\mathcal{B}(G)$ satisfies Hansen-Vuki{\v{c}}evi{\'c} conjecture if $|G|=2p, 2p^2, 4p$, $4p^2$ and $p^n$. While computing Zagreb indices of $\mathcal{B}(G)$ we have computed $\deg_{\mathcal{B}(G)}(H)$ for all $H \in L(G)$ for the above mentioned groups. Using these information, in Section 4, we  also compute Randic Connectivity index, Atom-Bond Connectivity index, Geometric-Arithmetic index, Harmonic index and Sum-Connectivity index  of $\mathcal{B}(G)$.
 \section{Structures of $\mathcal{B}(G)$}
 In \cite{DEN-23}, Das et al. obtained the structures of $\mathcal{B}(G)$ when $G$ is a cyclic group of order $2p$ or $2p^2$, where $p$ is an odd prime. We list the following results that are useful in our computation. We write $mK_{1, r}$ to denote $m$ copies of the star $K_{1, r}$.
\begin{theorem}\label{structure_cyclic_2p_2p^2}
 \cite[Observation 1.1 (b)]{DEN-23}   Let $G$ be a cyclic group and $p$ be an odd prime. 
 \begin{enumerate}
     \item If $|G|=2p$ then $\mathcal{B}(G) = K_2 \sqcup K_{1, 3}\sqcup K_{1, p^2 - 1} \sqcup K_{1, 3p^2 - 3}$. 
     \item  If $|G|=2p^2$ then $\mathcal{B}(G) = K_2 \sqcup K_{1, 3} \sqcup K_{1, p^2 - 1} \sqcup K_{1, 3p^2 - 3} \sqcup K_{1, p^4 - p^2} \sqcup K_{1, 3p^4 - 3p^2}$.
 \end{enumerate}
\end{theorem}
Das et al. \cite{DEN-23} also obtained the structures of $\mathcal{B}(G)$ when $G$ is the dihedral group of order $2p$ or $2p^2$ and the dicyclic group of order $4p$ or $4p^2$ (where $p$ is any prime).
\begin{theorem}\label{structure_of_D_2p}
\cite[Theorem 4.1-- 4.2]{DEN-23} Let $D_{2n}=\langle a, b: a^n=b^2=1, bab=a^{-1} \rangle$ be the dihedral group of order $2n$ ($n \geq 3$) and $p$ be a prime. Then 
	\[
	\mathcal{B}(D_{2n})= \begin{cases}
		K_2 \sqcup pK_{1, 3} \sqcup K_{1, p^2-1} \sqcup K_{1, 3p(p-1)}, & \text{ if } n = p\\
		K_2 \sqcup p^2K_{1, 3} \sqcup K_{1, p^2-1} \sqcup K_{1, p^4-p^2} \sqcup pK_{1, 3p(p-1)} \sqcup K_{1, 3p^2(p^2-p)}, & \text{ if } n = p^2.
	\end{cases}
	\] 
 \end{theorem}
 \begin{theorem}\label{structure_of_Q_4p}
 \cite[Theorem 4.3]{DEN-23}	Let $Q_{4p} =  \langle a, b : a^{2p} = 1, b^2 = a^p, bab^{-1} = a^{-1} \rangle$ be the dicyclic group of order $4p$, where $p$ is a prime. Then
 	\[
 	\mathcal{B}(Q_{4p})=\begin{cases}
 		K_2 \sqcup K_{1, 3} \sqcup 3K_{1, 12} \sqcup K_{1, 24}, & \text{ when } p=2 \\
 		K_2 \sqcup K_{1, 3} \sqcup pK_{1, 12} \sqcup K_{1, p^2-1} \sqcup K_{1, 3p^2-3} \sqcup K_{1, 12p^2-12p}, & \text{ when } p \geq 3.
 	\end{cases}
 	\]
 \end{theorem}
 \begin{theorem}\label{structure_of_Q_4p^2}
 \cite[Theorem 4.4]{DEN-23}	Let $Q_{4p^2} =  \langle a, b : a^{2p^2} = 1, b^2 = a^{p^2}, bab^{-1} = a^{-1} \rangle$ be the dicyclic group of order $4p^2$, where $p$ is a prime. Then
 	\[
 	\mathcal{B}(Q_{4p^2})=\begin{cases}
 		K_2 \sqcup K_{1, 3} \sqcup 5K_{1, 12} \sqcup 2K_{1, 24} \sqcup K_{1, 48} \sqcup K_{1, 96}, & \text{ when } p=2 \\
 		K_2 \sqcup K_{1, 3} \sqcup p^2K_{1, 12} \sqcup K_{1, p^2-1} \sqcup K_{1, 3p^2-3} \sqcup K_{1, 3p^4-3p^2} \\ \qquad \qquad \sqcup (p-1)K_{1, 12p^2-12p} \sqcup K_{1, 13p^4-12p^3+11p^2-12p}, & \text{ when } p \geq 3.
 	\end{cases}
 	\]
 \end{theorem}
In the following result we obtain the structures of $\mathcal{B}(G)$ when $G$ is  a cyclic group of order  $p^n$, $4p$ and  $4p^2$ for  any prime $p$. We write $\Nbd_{\mathcal{G}}(x)$ to denote the neighborhood of $x \in V(\mathcal{G})$, i. e., $\Nbd_{\mathcal{G}}(x) = \{y \in V(\mathcal{G}): y \text{ is adjacent to } x \in V(\mathcal{G})\}$. Also, $\mathcal{G}[S]$ denotes the subgraph of $\mathcal{G}$ induced by a subset $S$ of $V(\mathcal{G})$.
\begin{theorem}\label{structure_cyclic_p^n}
	If $G$ is a cyclic group of order $p^n$, where $p$ is a prime and $n \geq 1$, then $\mathcal{B}(G)=K_2 \sqcup K_{1, p^2-1} \sqcup K_{1, p^2(p^2-1)} \sqcup \ldots \sqcup K_{1, p^{2n-2}(p^2-1)}$.
\end{theorem}
\begin{proof}
Let $G = \langle a \rangle$. Then	 $V(\mathcal{B}(G))=G \times G \sqcup \{\{1\}, \langle a^{p^{n-1}} \rangle, \langle a^{p^{n-2}} \rangle, \langle a^{p^{n-3}} \rangle, \ldots, \langle a^{p} \rangle, \langle a \rangle\}$ as there are exactly $n+1$ subgroups of $G$ and they are cyclic. Here $|\langle a^{p^{n-1}} \rangle|=p, |\langle a^{p^{n-2}} \rangle|=p^2, |\langle a^{p^{n-3}} \rangle|=p^3, \ldots, |\langle a^{p} \rangle|=p^{n-1}$ \, and \, $|\langle a \rangle|=p^n$. As such, \, $\mathcal{B}(G)[\{\langle a^{p^{n-1}} \rangle\} \sqcup $ $\Nbd_{\mathcal{B}(G)}(\langle a^{p^{n-1}} \rangle)] = K_{1, p^2-p^{2-2}}=K_{1, p^2-1}$ and $\mathcal{B}(G)[\{\langle a^{p^{n-2}} \rangle\} \sqcup \Nbd_{\mathcal{B}(G)}(\langle a^{p^{n-2}} \rangle)] = K_{1, p^4-p^{4-2}}=K_{1, p^4-p^2}$. Similarly, $\mathcal{B}(G)[\{\langle a^{p^{n-3}} \rangle\} \sqcup \Nbd_{\mathcal{B}(G)}(\langle a^{p^{n-3}} \rangle)] = K_{1, p^6-p^{6-2}}=K_{1, p^6-p^4}$, $\mathcal{B}(G)[\{\langle a^{p^{n-4}} \rangle\} \sqcup \Nbd_{\mathcal{B}(G)}(\langle a^{p^{n-4}} \rangle)] = K_{1, p^8-p^{8-2}}=K_{1, p^8-p^6}$, \ldots, $\mathcal{B}(G)[\{\langle a^{p} \rangle\} \sqcup \Nbd_{\mathcal{B}(G)}(\langle a^{p} \rangle)] = K_{1, p^{2n-2}-p^{2n-2-2}}=K_{1, p^{2n-2}-p^{2n-4}}$ and $\mathcal{B}(G)[\{\langle a \rangle\} \sqcup \Nbd_{\mathcal{B}(G)}(\langle a \rangle)] = K_{1, p^{2n}-p^{2n-2}}$. Therefore,
	\begin{align*}
		\mathcal{B}(G)&= \underset{H \in L(G)}{\sqcup} \mathcal{B}(G)[\{H\} \sqcup \Nbd_{\mathcal{B}(G)}(H)] 
		= K_2 \sqcup K_{1, p^2-1} \sqcup K_{1, p^4-p^2} \sqcup \ldots \sqcup K_{1, p^{2n}-p^{2n-2}}.
	\end{align*}
~This completes the proof.
\end{proof}

\begin{theorem}\label{structure_cyclic_4p_4p^2}
	Let $G$ be a cyclic group and $p$ be an odd prime. 
	\begin{enumerate}
		\item If $|G|=4p$ then $\mathcal{B}(G) = K_2 \sqcup K_{1, 3} \sqcup K_{1, 12} \sqcup K_{1, p^2 - 1} \sqcup K_{1, 3p^2 - 3} \sqcup K_{1, 12p^2-12}$. 
		\item  If $|G|=4p^2$ then $\mathcal{B}(G) = K_2 \sqcup K_{1, 3} \sqcup K_{1, 12} \sqcup K_{1, p^2 - 1}  \sqcup K_{1, p^4 - p^2} \sqcup K_{1, 3p^2 - 3}$ $ \sqcup K_{1, 3p^4 - 3p^2} \sqcup K_{1, 12p^2-12} \sqcup K_{1, 12p^4-12p^2}$.
	\end{enumerate}
\end{theorem}
\begin{proof}
For any prime $p$, if $|G|=p^2$ then by Theorem \ref{structure_cyclic_p^n} we have
\begin{center}
	$\mathcal{B}(G)=K_2 \sqcup  K_{1, p^2 - 1}$ $ \sqcup K_{1, p^4 - p^2}$.
\end{center}

(a)  If $|G|=4p$ and $G = \langle a \rangle$ then $V(\mathcal{B}(G)) = G \times G \sqcup \{\{1\}, \langle a^{2p} \rangle, \langle a^p \rangle, \langle a^4 \rangle, \langle a^2 \rangle, \langle a \rangle\}$. Since  $|\langle a^{2p} \rangle| = 2$, $|\langle a^p \rangle| =4, |\langle a^4 \rangle|=p$ and $|\langle a^2 \rangle|=2p$ we have \, $\mathcal{B}(G)[\{\langle a^{2p} \rangle\} \sqcup \Nbd_{\mathcal{B}(G)}(\langle a^{2p} \rangle)] \, = \,  K_{1, 3}$, \quad $\mathcal{B}(G)[\{\langle a^{p} \rangle\} \sqcup \Nbd_{\mathcal{B}(G)}(\langle a^{p} \rangle)] \, = \, K_{1, 12}$, \quad  $\mathcal{B}(G)[\{\langle a^{4} \rangle\}$ $\sqcup$   $\Nbd_{\mathcal{B}(G)}(\langle a^{4} \rangle)] = K_{1, p^2-1}$ and $\mathcal{B}(G)[\{\langle a^2 \rangle\} \sqcup \Nbd_{\mathcal{B}(G)}(\langle a^2 \rangle)] = K_{1, 3p^2 - 3}$. Also, $\mathcal{B}(G)[\{\langle a \rangle\}$ $\sqcup \Nbd_{\mathcal{B}(G)}(\langle a \rangle)] = K_{1, 12p^2 - 12}$ noting that $|\Nbd_{\mathcal{B}(G)}(\langle a \rangle)| = 16p^2 - (1+3+12+p^2-1+3p^2-3) = 12p^2 - 12$. Thus, $\mathcal{B}(G) = K_2 \sqcup K_{1, 3} \sqcup K_{1, 12} \sqcup K_{1, p^2 - 1} \sqcup K_{1, 3p^2 - 3} \sqcup K_{1, 12p^2-12}$.

(b) If $|G|=4p^2$ and $G = \langle a \rangle$ then $V(\mathcal{B}(G)) = G \times G \sqcup \{\{1\}, \langle a^{2p^2} \rangle, \langle a^{p^2} \rangle, \langle a^{4p} \rangle, \langle a^4 \rangle, \langle a^{2p} \rangle,$ $\langle a^2 \rangle$, $\langle a^p \rangle, \langle a \rangle\}$. Since  $|\langle a^{2p^2} \rangle| = 2, |\langle a^{p^2} \rangle| = 4, |\langle a^{4p} \rangle| = p, |\langle a^{4} \rangle| = p^2, |\langle a^{2p} \rangle| = 2p$, $|\langle a^{2} \rangle| = 2p^2$ \, and \,  $|\langle a^p \rangle| =4p$ \, we have \, $\mathcal{B}(G)[\{\langle a^{2p^2} \rangle\} \sqcup \Nbd_{\mathcal{B}(G)}(\langle a^{2p^2} \rangle)] = K_{1, 3}$,  $\mathcal{B}(G)[\{\langle a^{p^2} \rangle\} \sqcup \Nbd_{\mathcal{B}(G)}(\langle a^{p^2} \rangle)] = K_{1, 12}, \quad  \mathcal{B}(G)[\{\langle a^{4p} \rangle\} \sqcup \Nbd_{\mathcal{B}(G)}(\langle a^{4p} \rangle)] = K_{1, p^2-1}$, \quad $\mathcal{B}(G)[\{\langle a^{4} \rangle\}$ $ \sqcup \Nbd_{\mathcal{B}(G)}(\langle a^{4} \rangle)] = K_{1, p^4-p^2}, \mathcal{B}(G)[\{\langle a^{2p} \rangle\} \sqcup \Nbd_{\mathcal{B}(G)}(\langle a^{2p} \rangle)] = K_{1, 3p^2-3}$, $\mathcal{B}(G)[\{\langle a^{2} \rangle\} \sqcup \Nbd_{\mathcal{B}(G)}(\langle a^{2} \rangle)] = K_{1, 3p^4-3p^2}$ and $\mathcal{B}(G)[\{\langle a^p \rangle\} \sqcup \Nbd_{\mathcal{B}(G)}(\langle a^p \rangle)] = K_{1, 12p^2 - 12}$. Also, \, $\mathcal{B}(G)[\{\langle a \rangle\} \sqcup \Nbd_{\mathcal{B}(G)}(\langle a \rangle)] = K_{1, 12p^4 - 12p^2}$ \, noting that 
\small{
\begin{align*}
 |\Nbd_{\mathcal{B}(G)}(\langle a \rangle)| &= 16p^4 - (1+3+12+p^2-1+p^4-p^2+3p^2-3+3p^4-3p^2+12p^2-12)\\
 & = 12p^4 - 12p^2.
\end{align*}} 
Thus, $\mathcal{B}(G) = K_2 \sqcup K_{1, 3} \sqcup K_{1, 12} \sqcup K_{1, p^2 - 1}  \sqcup K_{1, p^4 - p^2} \sqcup K_{1, 3p^2 - 3} \sqcup K_{1, 3p^4 - 3p^2} \sqcup K_{1, 12p^2-12} \sqcup K_{1, 12p^4-12p^2}$.
\end{proof}

We conclude this section noting that   $\mathcal{B}(G)$ is a union of   star graphs and number of components in $\mathcal{B}(G)$ is $|L(G)|$ if the number of generators of $G$ is at most two. For each subgroup $H$ of $G$ we have $\mathcal{B}(G)[\{H\} \sqcup \Nbd(H)] = K_{1, m}$, for some $m$.  

\section{Zagreb indices of $\mathcal{B}(G)$}
In this section, we obtain expressions for first and second   Zagreb indices of $\mathcal{B}(G)$ and  obtain a condition such that $\mathcal{B}(G)$ satisfy Hansen-Vuki{\v{c}}evi{\'c} conjecture \cite{hansen2007comparing}. Using this we show that $\mathcal{B}(G)$ satisfy Hansen-Vuki{\v{c}}evi{\'c} conjecture for the groups considered in Section 2. 
We begin with the following result.
\begin{theorem}\label{Zagreb_indices_B(G)}
	Let $G$ be a finite group. Then the Zagreb indices of $\mathcal{B}(G)$ are given by
	\[
	M_1(\mathcal{B}(G))=|G|^2+|G|^4 \sum\limits_{H \in L(G)}\left({\Pr}_H(G)\right)^2 \quad \text{and} \quad M_2(\mathcal{B}(G))= M_1(\mathcal{B}(G))-|G|^2.
	\]
\end{theorem}
\begin{proof}
	For $(a, b) \in G \times G$ and $H \in L(G)$, by definition of $\mathcal{B}(G)$   and \eqref{deg(H in L(G))}, we have 
	\[
	\deg_{\mathcal{B}(G)}((a, b))=1 \quad \text{and} \quad \deg_{\mathcal{B}(G)}(H)=|G|^2 {\Pr}_H(G).
	\]
	Therefore, 
	\begin{align*}
		M_1(\mathcal{B}(G))&=\sum_{(a, b) \in G \times G}\left(\deg_{\mathcal{B}(G)}((a, b))\right)^2+\sum_{H \in L(G)}\left(\deg_{\mathcal{B}(G)}(H)\right)^2\\
		&= \sum_{(a, b) \in G \times G} 1 + \sum_{H \in L(G)}\left(|G|^2{\Pr}_H(G)\right)^2 
		= |G|^2+|G|^4 \sum_{H \in L(G)}\left({\Pr}_H(G)\right)^2
	\end{align*}
	and
	\begin{align*}
		M_2(\mathcal{B}(G))&=\sum_{(a, b)H \in e(\mathcal{B}(G))}\deg_{\mathcal{B}(G)}((a, b))\deg_{\mathcal{B}(G)}(H) \\
		&= \sum_{(a, b)H \in e(\mathcal{B}(G))} \deg_{\mathcal{B}(G)}(H) 
		= \sum_{H \in L(G)} \left(\deg_{\mathcal{B}(G)}(H)\right)^2, 
	\end{align*}
since for each edge incident to $H$, $\deg_{\mathcal{B}(G)}(H)$ is added.
Therefore, 
\[
M_2(\mathcal{B}(G))= \sum_{H \in L(G)}\left(|G|^2{\Pr}_H(G)\right)^2 
	= |G|^4 \sum_{H \in L(G)}\left({\Pr}_H(G)\right)^2.
\]
	Hence, the result follows.
\end{proof} 

Now we derive a necessary and sufficient condition such that $\mathcal{B}(G)$ satisfy Hansen-Vuki{\v{c}}evi{\'c} conjecture.  The following result, which gives the size of $\mathcal{B}(G)$, is useful.
\begin{theorem}\label{deg_sum=num_of_edges}
	\cite[Theorem 3.2]{DEN-23} Let $G$ be a finite group and $e(\mathcal{B}(G))$ be the set of edges of the graph $\mathcal{B}(G)$. Then

	\[
	\sum_{(a, b)\in G\times G} \deg_{\mathcal{B}(G)}((a, b))= \sum_{H \in L(G)} \deg_{\mathcal{B}(G)}(H)= |G|^2=|e(\mathcal{B}(G))|.
	\]
\end{theorem}   

\begin{theorem}\label{condition_for_satisfying_conjecture} 
	Let $G$ be a finite group. Then \quad $\frac{M_{2}(\mathcal{B}(G))}{\vert e(\mathcal{B}(G))\vert} \geq \frac{M_{1}(\mathcal{B}(G))}{\vert V(\mathcal{B}(G)) \vert} $ if and only if  $|L(G)|P-1 \geq 0$, where $P=\sum\limits_{H \in L(G)}(\Pr_H(G))^2$.
\end{theorem}
\begin{proof}
	We have $|V(\mathcal{B}(G))|=|G|^2+|L(G)|$ and from Theorem \ref{deg_sum=num_of_edges}, $|e(\mathcal{B}(G))|=|G|^2$. Using Theorem \ref{Zagreb_indices_B(G)}, we get 
	
	\begin{align*}
		\frac{M_{2}(\mathcal{B}(G))}{\vert e(\mathcal{B}(G))\vert} - \frac{M_{1}(\mathcal{B}(G))}{\vert V(\mathcal{B}(G)) \vert}&= \frac{M_1(\mathcal{B}(G))-|G|^2}{|G|^2}-\frac{M_{1}(\mathcal{B}(G))}{|G|^2+|L(G)|} \\
		&= \frac{|L(G)|(M_1(\mathcal{B}(G)-|G|^2)-|G|^4}{|G|^2(|G|^2+|L(G)|)} \\
		&= \frac{|L(G)||G|^4  \sum\limits_{H \in L(G)}({\Pr}_H(G))^2 - |G|^4}{|G|^2(|G|^2+|L(G)|)} = \frac{|G|^4(|L(G)|P-1)}{|G|^2(|G|^2+|L(G)|)},
	\end{align*}
	where \quad $P=\sum\limits_{H \in L(G)}(\Pr_H(G))^2$. Hence, \quad  $\frac{M_{2}(\mathcal{B}(G))}{\vert e(\mathcal{B}(G))\vert} \geq \frac{M_{1}(\mathcal{B}(G))}{\vert V(\mathcal{B}(G)) \vert} $ \quad if and only if $|L(G)|P-1 \geq 0$.
\end{proof}

In the rest part of this section we shall consider  $\mathcal{B}(G)$ when $G$ is a cyclic group, dihedral group and dicyclic group  of certain orders. 
\begin{theorem}\label{cyclic case}
Let $G$ be a cyclic group and $p$ be an odd prime. 
	\begin{enumerate}
		\item If $|G|=2p$ then $M_1(\mathcal{B}(G))=10p^4-16p^2+20$ and $M_2(\mathcal{B}(G))=10p^4-20p^2+20$.
		\item If $|G|=2p^2$ then $M_1(\mathcal{B}(G))=10p^8-20p^6+24p^4-20p^2+20$ and $M_2(\mathcal{B}(G))=10p^8-20p^6+20p^4-20p^2+20$.
		\item If $|G|=4p$ then $M_1(\mathcal{B}(G))=154p^4-292p^2+308$ and $M_2(\mathcal{B}(G))=154p^4-308p^2+308$.
		\item If $|G|=4p^2$ then $M_1(\mathcal{B}(G))=154p^8-308p^6+324p^4-308p^2+308$ and $M_2(\mathcal{B}(G))=154p^8-308p^6+308p^4-308p^2+308$.
	\end{enumerate}
	Further, $\frac{M_{2}(\mathcal{B}(G))}{\vert e(\mathcal{B}(G))\vert} \geq \frac{M_{1}(\mathcal{B}(G))}{\vert V(\mathcal{B}(G)) \vert}$.
\end{theorem}
\begin{proof}
	(a) If $|G|=2p$ then, by Theorem \ref{structure_cyclic_2p_2p^2}(a), we have $\mathcal{B}(G)=K_2 \sqcup K_{1, 3} \sqcup K_{1, p^2-1} \sqcup K_{1, 3p^2-3}$. As such $|L(G)| = 4$ and $P=\sum\limits_{H\in L(G)}(\Pr_H(G))^2=\sum\limits_{i=1}^{4}\left(\Pr_{H_i}(G)\right)^2$. By \eqref{deg(H in L(G))}, we get
		\begin{align*}
			|G|^4\sum_{i=1}^{4}\left({\Pr}_{H_i}(G)\right)^2 &=\sum_{i=1}^{4}\left(\deg_{\mathcal{B}(G)}(H_i)\right)^2 \\
			&= 1^2+3^2+(p^2-1)^2+(3p^2-3)^2 
			= 10p^4-20p^2+20.
		\end{align*} 
  Therefore, by Theorem \ref{Zagreb_indices_B(G)}, we have
\[
M_1(\mathcal{B}(G))=4p^2+ 10p^4-20p^2+20 
           = 10p^4-16p^2+20
\]
and
\[
M_2(\mathcal{B}(G))= 10p^4-16p^2+20-4p^2
                    = 10p^4-20p^2+20.
\]
		Further, 
		\[
		|L(G)|P-1= \frac{4(10p^4-20p^2+20)}{16p^4}-1= \frac{2p^2(3p^2-10)+20}{4p^4}. 
		\]
		If $p >2$ then $|L(G)|P-1 > \frac{9}{16}$. Thus,  $|L(G)|P-1 > 0$ and hence the result follows from Theorem \ref{condition_for_satisfying_conjecture}.
		
		 (b) If $|G|=2p^2$ then, by Theorem \ref{structure_cyclic_2p_2p^2}(b), we have $\mathcal{B}(G)=K_2 \sqcup K_{1, 3} \sqcup K_{1, p^2 - 1} \sqcup K_{1, 3p^2 - 3} \sqcup K_{1, p^4 - p^2} \sqcup K_{1, 3p^4 - 3p^2}$. As such  $|L(G)| = 6$ and $P=\sum\limits_{H\in L(G)}(\Pr_H(G))^2=\sum\limits_{i=1}^{6}\left(\Pr_{H_i}(G)\right)^2$. By \eqref{deg(H in L(G))}, we get
		\begin{align*}
			|G|^4\sum_{i=1}^{6}\left({\Pr}_{H_i}(G)\right)^2 &=\sum_{i=1}^{6}\left(\deg_{\mathcal{B}(G)}(H_i)\right)^2 \\
			&=1^2+3^2+(p^2-1)^2+(3p^2-3)^2+(p^4-p^2)^2+(3p^4-3p^2)^2  \\
			&= 10p^8-20p^6+20p^4-20p^2+20.
		\end{align*} 
  Therefore, by Theorem \ref{Zagreb_indices_B(G)}, we have
		\begin{align*}
			M_1(\mathcal{B}(G))&= 4p^4+10p^8-20p^6+20p^4-20p^2+20\\
                &= 10p^8-20p^6+24p^4-20p^2+20
		\end{align*}
		and
		\begin{align*}
			M_2(\mathcal{B}(G))&=10p^8-20p^6+24p^4-20p^2+20-4p^4 \\
                &= 10p^8-20p^6+20p^4-20p^2+20.
		\end{align*}
		Further, 
		\begin{align*}
			|L(G)|P-1&= \frac{6(10p^8-20p^6+20p^4-20p^2+20)}{16p^8}-1 \\
			&= \frac{p^6(22p^2-60)+60p^2(p^2-1)+60}{8p^8}. 
		\end{align*}
		We have $22p^2-60>0$ and $p^2-1>0$ for all $p > 2$. Thus,  $|L(G)|P-1 > 0$ and hence the result follows from Theorem \ref{condition_for_satisfying_conjecture}.

  (c) If $|G|=4p$ then, by Theorem \ref{structure_cyclic_4p_4p^2}(a), we have $\mathcal{B}(G)=K_2 \sqcup K_{1, 3} \sqcup K_{1, 12} \sqcup K_{1, p^2 - 1} \sqcup K_{1, 3p^2 - 3} \sqcup K_{1, 12p^2-12}$.  As such $|L(G)| = 6$ and $P=\sum\limits_{H\in L(G)}(\Pr_H(G))^2=\sum\limits_{i=1}^{6}\left(\Pr_{H_i}(G)\right)^2$. By \eqref{deg(H in L(G))}, we get
		\begin{align*}
			|G|^4\sum_{i=1}^{6}\left({\Pr}_{H_i}(G)\right)^2 &=\sum_{i=1}^{6}\left(\deg_{\mathcal{B}(G)}(H_i)\right)^2 \\
			&= 1^2+3^2+12^2+(p^2-1)^2+(3p^2-3)^2+(12p^2-12)^2 \\
			&= 154p^4-308p^2+308.
		\end{align*}
  Therefore, by Theorem \ref{Zagreb_indices_B(G)}, we have
		\begin{align*}
			M_1(\mathcal{B}(G))&=16p^2+154p^4-308p^2+308
			= 154p^4-292p^2+308 
		\end{align*}
		and
		\begin{align*}
			M_2(\mathcal{B}(G))&=154p^4-292p^2+308-16p^2
			= 154p^4-308p^2+308.
		\end{align*}
  Further,
        \[
		|L(G)|P-1= \frac{6(154p^4-308p^2+308)}{256p^4}-1= \frac{p^2(334p^2-924)+924}{128p^4}. 
		\]
		For all $p > 2$ we have $334p^2-924>0$. Thus,  $|L(G)|P-1 > 0$ and hence the result follows from Theorem \ref{condition_for_satisfying_conjecture}.

   (d) If $|G|=4p^2$ then, by Theorem \ref{structure_cyclic_4p_4p^2}(b), we have $\mathcal{B}(G)= K_2 \sqcup K_{1, 3} \sqcup K_{1, 12} \sqcup K_{1, p^2 - 1}  \sqcup K_{1, p^4 - p^2} \sqcup K_{1, 3p^2 - 3} \sqcup K_{1, 3p^4 - 3p^2} \sqcup K_{1, 12p^2-12} \sqcup K_{1, 12p^4-12p^2}$. As such $|L(G)| = 9$ and $P=\sum\limits_{H\in L(G)}(\Pr_H(G))^2=\sum\limits_{i=1}^{9}\left(\Pr_{H_i}(G)\right)^2$. By \eqref{deg(H in L(G))}, we get
		\begin{align*}
			|G|^4\sum_{i=1}^{9}\left({\Pr}_{H_i}(G)\right)^2 &=\sum_{i=1}^{9}\left(\deg_{\mathcal{B}(G)}(H_i)\right)^2 \\
			&=1^2+3^2+12^2+(p^2-1)^2+(p^4-p^2)^2+(3p^2 -3)^2 \\
			& \qquad \qquad \qquad  +(3p^4-3p^2)^2+(12p^2-12)^2 +(12p^4-12p^2)^2   \\
			&=154p^8-308p^6+308p^4-308p^2+308.
		\end{align*}
  Therefore, by Theorem \ref{Zagreb_indices_B(G)}, we have
		\begin{align*}
			M_1(\mathcal{B}(G))&=16p^4+154p^8-308p^6+308p^4-308p^2+308 \\
			&= 154p^8-308p^6+324p^4-308p^2+308
		\end{align*}
		and
		\begin{align*}
			M_2(\mathcal{B}(G))&=154p^8-308p^6+324p^4-308p^2+308-16p^4\\
			&= 154p^8-308p^6+308p^4-308p^2+308.
		\end{align*}
		Further, 
		\begin{align*}
			|L(G)|P-1&= \frac{9(154p^8-308p^6+308p^4-308p^2+308)}{256p^8}-1 \\
			&= \frac{p^6(1130p^2-2772)+2772p^2(p^2-1)+2772}{256p^8}. 
		\end{align*}
		We have $1130p^2-2772>0$ and $p^2-1>0$ for all $p > 2$. Thus,  $|L(G)|P-1 > 0$ and hence the result follows from Theorem \ref{condition_for_satisfying_conjecture}.
\end{proof}
 
If $p = 2$ then the groups considered in Theorem \ref{cyclic case} becomes cyclic groups of order $2^n$ where $n = 2, 3, 4$. These groups, more generally cyclic groups of order $p^n$, are considered in our next result.

\begin{theorem}
    If $G$ is a cyclic group of order $p^n$, where $p$ is a prime and $n \geq 1$, then $M_1(\mathcal{B}(G))=\frac{p^{2n}(p^2+1)+p^{4n}(p^2-1)+2}{p^2+1}$ and $ M_2(\mathcal{B}(G))=\frac{p^{4n}(p^2-1)+2}{p^2+1}$. Further, $\frac{M_{2}(\mathcal{B}(G))}{\vert e(\mathcal{B}(G))\vert} \geq \frac{M_{1}(\mathcal{B}(G))}{\vert V(\mathcal{B}(G)) \vert}$.
\end{theorem}
\begin{proof}
    By Theorem \ref{structure_cyclic_p^n}, we have $\mathcal{B}(G)=K_2 \sqcup K_{1, p^2-1} \sqcup K_{1, p^2(p^2-1)} \sqcup \ldots \sqcup K_{1, p^{2n-2}(p^2-1)}$. As such $|L(G)| = n+1$ and $P=\sum\limits_{H\in L(G)}(\Pr_H(G))^2=\sum\limits_{i=1}^{n+1}\left(\Pr_{H_i}(G)\right)^2$. By \eqref{deg(H in L(G))}, we get
    \begin{align*}
		|G|^4\sum_{i=1}^{n+1}\left({\Pr}_{H_i}(G)\right)^2 &=\sum_{i=1}^{n+1}\left(\deg_{\mathcal{B}(G)}(H_i)\right)^2 \\
		&= 1^2+(p^2-1)^2+(p^4-p^2)^2+(p^6-p^4)^2+\ldots+(p^{2n}-p^{2n-2})^2 \\
	&= 1^2+(p^2-1)^2\left(1+p^4+p^8+p^{12}+\ldots+p^{4n-4}\right) \\
        &=1+(p^2-1)^2 \left(\frac{p^{4n}-1}{p^4-1}\right)=\frac{p^{4n}(p^2-1)+2}{p^2+1}.
	\end{align*}
 Therefore, by Theorem \ref{Zagreb_indices_B(G)}, we have
\begin{align*}
M_1(\mathcal{B}(G))& = p^{2n}+\frac{p^{4n}(p^2-1)+2}{p^2+1}=\frac{p^{2n}(p^2+1)+p^{4n}(p^2-1)+2}{p^2+1}
\end{align*}
and
\begin{align*}
M_2(\mathcal{B}(G)) &= \frac{p^{2n}(p^2+1)+p^{4n}(p^2-1)+2}{p^2+1} - p^{2n}=\frac{p^{4n}(p^2-1)+2}{p^2+1}.
\end{align*}
Further, 
	\begin{align*}
		|L(G)|P-1= (n+1)\frac{p^{4n}(p^2-1)+2}{p^{4n}(p^2+1)}-1 &= \frac{np^{4n+2}-np^{4n}+2n-p^{4n}+2-p^{4n}}{p^{4n}(p^2+1)} \\
		&= \frac{p^{4n}(np^2-n-2)+2n+2}{p^{4n}(p^2+1)}. 
	\end{align*}
	For all $p \geq 2$ and $n \geq 1$, we have $n(p^2-1) >2$. Thus,  $|L(G)|P-1 > 0$ and hence the result follows from Theorem \ref{condition_for_satisfying_conjecture}.
\end{proof}

\begin{theorem}
	If $G=D_{2p}$, then $M_1(\mathcal{B}(G))=10p^4-18p^3 +11p^2+9p+2$ and $ M_2(\mathcal{B}(G))=10p^4-18p^3 +7p^2+9p+2$. Further, $\frac{M_{2}(\mathcal{B}(G))}{\vert e(\mathcal{B}(G))\vert} \geq \frac{M_{1}(\mathcal{B}(G))}{\vert V(\mathcal{B}(G)) \vert}$.
\end{theorem}
\begin{proof}
	By Theorem \ref{structure_of_D_2p}, we have $\mathcal{B}(G)=K_2 \sqcup pK_{1, 3} \sqcup K_{1, p^2-1} \sqcup K_{1, 3p(p-1)}$. As such $|L(D_{2p})| = p+3$ and $P=\sum\limits_{H\in L(G)}(\Pr_H(G))^2=\sum\limits_{i=1}^{p+3}\left(\Pr_{H_i}(G)\right)^2$. By \eqref{deg(H in L(G))}, we get
	\begin{align*}
		|G|^4\sum_{i=1}^{p+3}\left({\Pr}_{H_i}(G)\right)^2 &=\sum_{i=1}^{p+3}\left(\deg_{\mathcal{B}(G)}(H_i)\right)^2 \\
		&= 1^2+(\underbrace{3^2+\ldots+3^2}_{p\text{-times}})+(p^2-1)^2+(3p^2-3p)^2 \\
		&= 10p^4-18p^3 +7p^2+9p+2.
	\end{align*}
Therefore, by Theorem \ref{Zagreb_indices_B(G)}, we have
\begin{align*}
M_1(\mathcal{B}(G))& = 4p^2 + 10p^4-18p^3 +7p^2+9p+2\\
&= 10p^4-18p^3 +11p^2+9p+2
\end{align*}
	and
\begin{align*}
M_2(\mathcal{B}(G)) &= 10p^4-18p^3 +11p^2+9p+2 - 4p^2\\
&= 10p^4-18p^3 +7p^2+9p+2.
\end{align*}
Further, 
	\begin{align*}
		|L(G)|P-1&= \frac{(p+3)(10p^4-18p^3 +7p^2+9p+2)}{16p^4}-1 \\
		&= \frac{p^3(10p^2-4p-47)+30p^2+29p+6}{16p^4}. 
	\end{align*}
	If $p = 2$ then $|L(G)|P-1 = \frac{1}{4}$.  If $p \geq 3$ then $|L(G)|P-1 \geq \frac{25}{27}$. Thus,  $|L(G)|P-1 > 0$ and hence the result follows from Theorem \ref{condition_for_satisfying_conjecture}.
\end{proof}
\begin{theorem}
	If $G=D_{2p^2}$, then $M_1(\mathcal{B}(G))=10p^8-18p^7 +7p^6+9p^5-12p^4+9p^3+7p^2+2$ and $ M_2(\mathcal{B}(G))=10p^8-18p^7 +7p^6+9p^5-16p^4+9p^3+7p^2+2$. Further, $\frac{M_{2}(\mathcal{B}(G))}{\vert e(\mathcal{B}(G))\vert} \geq \frac{M_{1}(\mathcal{B}(G))}{\vert V(\mathcal{B}(G)) \vert}$.
\end{theorem}
\begin{proof}
	By Theorem \ref{structure_of_D_2p}, we have $\mathcal{B}(G)=K_2 \sqcup p^2K_{1, 3} \sqcup K_{1, p^2-1} \sqcup K_{1, p^4-p^2} \sqcup pK_{1, 3p(p-1)} \sqcup K_{1, 3p^2(p^2-p)}$. As such $|L(D_{2p^2})| = p^2+p+4$ and $P=\sum\limits_{H\in L(G)}(\Pr_H(G))^2=\sum\limits_{i=1}^{p^2+p+4}(\Pr_{H_i}(G))^2$. 
 
  By  \eqref{deg(H in L(G))}, we get
	\begin{align*}
		|G|^4\sum_{i=1}^{p^2+p+4}({\Pr}_{H_i}(G))^2 &=\sum_{i=1}^{p^2+p+4}\left(\deg_{\mathcal{B}(G)}(H_i)\right)^2 \\
		&=1^2+(\underbrace{3^2+\ldots+3^2}_{p^2\text{-times}})+(p^2-1)^2+(p^4-p^2)^2 \\
		&\quad \quad \quad +(\underbrace{(3p^2-3p)^2+\ldots+(3p^2-3p)^2}_{p\text{-times}})+(3p^4-3p^3)^2 \\
		&= 10p^8-18p^7 +7p^6+9p^5-16p^4+9p^3+7p^2+2.
	\end{align*}
 Therefore, by Theorem \ref{Zagreb_indices_B(G)}, we have
	\begin{align*}
		M_1(\mathcal{B}(G))&=4p^4+10p^8-18p^7 +7p^6+9p^5-16p^4+9p^3+7p^2+2 \\
		&= 10p^8-18p^7 +7p^6+9p^5-12p^4+9p^3+7p^2+2 
	\end{align*}
	and     
	\begin{align*}
		M_2(\mathcal{B}(G))&=10p^8-18p^7 +7p^6+9p^5-12p^4+9p^3+7p^2+2-4p^4\\
		&=10p^8-18p^7 +7p^6+9p^5-16p^4+9p^3+7p^2+2.
	\end{align*}
	Further, 
	\begin{align*}
		|L(G)|P-1&= \frac{(p^2+p+4)(10p^8-18p^7 +7p^6+9p^5-16p^4+9p^3+7p^2+2)}{16p^8}-1 \\
		&=\frac{1}{16p^8}\bigg(p^9(6p-8)+p^7(4p^3+13p-56)+21p^6+p^4(29p-48)\\
		&\quad \quad \quad \quad +43p^3+30p^2+2p+8\bigg).
	\end{align*}
	For all $p \geq 2$, we have  $6p-8>0$, $4p^3+13p-56>0$ and $29p-48>0$. Therefore, $|L(G)|P-1>0$ and so the result follows from Theorem \ref{condition_for_satisfying_conjecture}. 
\end{proof}
\begin{theorem}
Consider the group $G=Q_{4p}$, where $p$ is a prime.
    \begin{enumerate}
        \item If $p=2$ then $M_1(\mathcal{B}(G))=1082$ and $M_2(\mathcal{B}(G))=1018$.
        \item If $p \geq 3$ then $M_1(\mathcal{B}(G))=154p^4-288p^3+140p^2+144p+20$ and $M_2(\mathcal{B}(G))=154p^4-288p^3+124p^2+144p+20$.
    \end{enumerate}
    Further, $\frac{M_{2}(\mathcal{B}(G))}{\vert e(\mathcal{B}(G))\vert} \geq \frac{M_{1}(\mathcal{B}(G))}{\vert V(\mathcal{B}(G)) \vert}$.
\end{theorem}
\begin{proof}
   (a) If $p=2$ then, by Theorem \ref{structure_of_Q_4p}, we have $\mathcal{B}(G)=K_2 \sqcup K_{1, 3} \sqcup 3K_{1, 12} \sqcup K_{1, 24}$. As such $|L(G)| = 6$ and $P=\sum\limits_{H\in L(G)}(\Pr_H(G))^2=\sum\limits_{i=1}^{6}\left(\Pr_{H_i}(G)\right)^2$. By \eqref{deg(H in L(G))}, we get
	\begin{align*}
		|G|^4\sum_{i=1}^{6}\left({\Pr}_{H_i}(G)\right)^2 &=\sum_{i=1}^{6}\left(\deg_{\mathcal{B}(G)}(H_i)\right)^2 \\
        &=1^2+3^2+12^2+12^2+12^2+24^2 
        = 1018.
	\end{align*}
 Therefore, by Theorem \ref{Zagreb_indices_B(G)}, we have
       \begin{align*}
           M_1(\mathcal{B}(G))=64+1018=1082   \text{ and }
           M_2(\mathcal{B}(G))=1082-64=1018.
       \end{align*}
       Further, $|L(G)|P-1= \frac{2012}{4096} >0$ and hence the result follows from Theorem \ref{condition_for_satisfying_conjecture}.

   (b) If $p \geq 3$ then, by Theorem \ref{structure_of_Q_4p}, we have $\mathcal{B}(G)=K_2 \sqcup K_{1, 3} \sqcup pK_{1, 12} \sqcup K_{1, p^2-1} \sqcup K_{1, 3p^2-3} \sqcup K_{1, 12p^2-12p}$. As such $|L(G)| = p+5$ and $P=\sum\limits_{H\in L(G)}(\Pr_H(G))^2=\sum\limits_{i=1}^{p+5}\left(\Pr_{H_i}(G)\right)^2$. By \eqref{deg(H in L(G))}, we get
	\begin{align*}
		|G|^4\sum_{i=1}^{p+5}&\left({\Pr}_{H_i}(G)\right)^2 =\sum_{i=1}^{p+5}\left(\deg_{\mathcal{B}(G)}(H_i)\right)^2 \\
		&= 1^2+3^2+(\underbrace{12^2+\ldots+12^2}_{p-\text{times}})+(p^2-1)^2+(3p^2-3)^2+(12p^2-12p)^2 \\
		&= 154p^4-288p^3+124p^2+144p+20.
	\end{align*}
 Therefore, by Theorem \ref{Zagreb_indices_B(G)}, we have
 \begin{align*}
     M_1(\mathcal{B}(G))&=16p^2+ 154p^4-288p^3+124p^2+144p+20\\
		&= 154p^4-288p^3+140p^2+144p+20 
	\end{align*}
	and
	\begin{align*}
		M_2(\mathcal{B}(G))&= 154p^4-288p^3+140p^2+144p+20 -16p^2\\
		&= 154p^4-288p^3+124p^2+144p+20.
	\end{align*}
	Further, 
	\begin{align*}
		|L(G)|P-1&= \frac{(p+5)(154p^4-288p^3+124p^2+144p+20)}{256p^4}-1 \\
            &= \frac{p^3(154p^2+226p-1316)+764p^2+740p+100}{256p^4}. 
	\end{align*}
	For all $p > 2$ we have $154p^2+226p>1316$. Thus,  $|L(G)|P-1 > 0$ and hence the result follows from Theorem \ref{condition_for_satisfying_conjecture}.
\end{proof}

\begin{theorem}\label{dycyclic-2}
Consider the group $G=Q_{4p^2}$, where $p$ is a prime. 
    \begin{enumerate}
        \item If $p=2$ then $M_1(\mathcal{B}(G))=13658$ and $M_2(\mathcal{B}(G))=13402$.
        \item If $p \geq 3$ then $M_1(\mathcal{B}(G))=178p^8-312p^7+412p^6-432p^5+12p^4+168p^3+124p^2+20$ and $M_2(\mathcal{B}(G))=178p^8-312p^7+412p^6-432p^5-4p^4+168p^3+124p^2+20$.
    \end{enumerate}
    Further, $\frac{M_{2}(\mathcal{B}(G))}{\vert e(\mathcal{B}(G))\vert} \geq \frac{M_{1}(\mathcal{B}(G))}{\vert V(\mathcal{B}(G)) \vert}$.
\end{theorem}
\begin{proof}
   (a) If $p=2$ then, by Theorem \ref{structure_of_Q_4p^2}, we have $\mathcal{B}(G)=K_2 \sqcup K_{1, 3} \sqcup 5K_{1, 12} \sqcup 2K_{1, 24} \sqcup K_{1, 48} \sqcup K_{1, 96}$. As such $|L(G)| = 11$ and $P=\sum\limits_{H\in L(G)}(\Pr_H(G))^2=\sum\limits_{i=1}^{11}\left(\Pr_{H_i}(G)\right)^2$. By \eqref{deg(H in L(G))}, we get
	\begin{align*}
		|G|^4\sum_{i=1}^{11}\left({\Pr}_{H_i}(G)\right)^2 &=\sum_{i=1}^{11}\left(\deg_{\mathcal{B}(G)}(H_i)\right)^2 \\
        &=1^2+3^2+(\underbrace{12^2+\ldots+12^2}_{5-\text{times}})+24^2 +24^2+48^2+96^2\\
        &=13402.
	\end{align*}
 Therefore, by Theorem \ref{Zagreb_indices_B(G)}, we have
\[           
M_1(\mathcal{B}(G))=256+13402=13658  \text{ and }      
           M_2(\mathcal{B}(G))=13658-256=13402.
\]
       Further, $|L(G)|P-1= \frac{81886}{65536} >0$ and hence the result follows from Theorem \ref{condition_for_satisfying_conjecture}.
       
       (b) If $p \geq 3$ then, by Theorem \ref{structure_of_Q_4p^2}, we have $\mathcal{B}(G)=K_2 \sqcup K_{1, 3} \sqcup p^2K_{1, 12} \sqcup K_{1, p^2-1} \sqcup K_{1, 3p^2-3} \sqcup K_{1, 3p^4-3p^2}  \sqcup (p-1)K_{1, 12p^2-12p} \sqcup K_{1, 13p^4-12p^3+11p^2-12p}$. As such $|L(G)| =p^2+ p+5$ and $P=\sum\limits_{H\in L(G)}(\Pr_H(G))^2=\sum\limits_{i=1}^{p^2+p+5}\left(\Pr_{H_i}(G)\right)^2$. By \eqref{deg(H in L(G))}, we get
	\begin{align*}
		|G|^4\sum_{i=1}^{p^2+p+5}\left({\Pr}_{H_i}(G)\right)^2 &=\sum_{i=1}^{p^2+p+5}\left(\deg_{\mathcal{B}(G)}(H_i)\right)^2 \\
		&=1^2+3^2+(\underbrace{12^2+\ldots+12^2}_{p^2-\text{times}})+(p^2-1)^2+(3p^2-3)^2\\
            & \qquad \quad +(3p^4-3p^2)^2 +(\underbrace{(12p^2-12p)^2+\ldots+(12p^2-12p)^2}_{(p-1)-\text{times}}) \\
            & \qquad \qquad \qquad \qquad \qquad +(13p^4-12p^3+11p^2-12p)^2  \\
		&=178p^8-312p^7+412p^6-432p^5-4p^4+168p^3+124p^2+20.
	\end{align*}
 Therefore, by Theorem \ref{Zagreb_indices_B(G)}, we have
 \begin{align*}
     M_1(\mathcal{B}(G))&=16p^4+178p^8-312p^7+412p^6-432p^5-4p^4+168p^3+124p^2+20 \\
		&= 178p^8-312p^7+412p^6-432p^5+12p^4+168p^3+124p^2+20
	\end{align*}
	and
	\begin{align*}
		M_2(\mathcal{B}(G))&=178p^8-312p^7+412p^6-432p^5+12p^4+168p^3+124p^2+20+16p^4\\
		&= 178p^8-312p^7+412p^6-432p^5-4p^4+168p^3+124p^2+20.
	\end{align*}
	Further, 
	\begin{align*}
		|L(G)|&P-1 \\
            & = \frac{(p^2+p+5)(178p^8-312p^7+412p^6-432p^5-4p^4+168p^3+124p^2+20)}{256p^8}-1 \\
            &= \frac{p^9(178p-134)+p^7(734p-1580)+p^5(1624p-1996)}{256p^8} \\
            & \qquad  \qquad \qquad \qquad \qquad \qquad \qquad \qquad \qquad +\frac{272p^4+964p^3+640p^2+20p+100}{256p^8} . 
	\end{align*}
	For all $p \geq 3$ we have $178p-134>0, 734p-1580 > 0$ and $1624p-1996>0$. Thus,  $|L(G)|P-1 > 0$ and hence the result follows from Theorem \ref{condition_for_satisfying_conjecture}.
\end{proof}

We conclude this section noting that $\mathcal{B}(G)$ satisfies Hansen-Vuki{\v{c}}evi{\'c} conjecture if $G$ is isomorphic to
\begin{enumerate}
    \item the cylic groups of order $2p, 2p^2, 4p$,  $4p^2$ and $p^n$ for any  prime $p$ and $n \geq 1$.
    \item the dihedral groups of order $2p$ and $2p^2$ for any prime $p$.
    \item the dicyclic groups of order $4p$ and $4p^2$ for any prime $p$.
\end{enumerate}
We have not found  any counterexample as yet. Therefore, the following problem is worth considering.
\begin{prob}
Is there any finite group $G$ such that $\mathcal{B}(G)$	does not satisfy Hansen-Vuki{\v{c}}evi{\'c} conjecture?
\end{prob}

\section{Other topological indices}
Several other degree-based topological indices were introduced soon after Zagreb indices. Randic Connectivity index, Atom-Bond Connectivity index, Geometric-Arithmetic index, Harmonic index, Sum-Connectivity index etc. are among the popular ones. A survey on these degree-based topological indices can be found in \cite{MNJ-FA-2020}. The Randic Connectivity index $R(\mathcal{G})$, Atom-Bond Connectivity index $\ABC(\mathcal{G})$, Geometric-Arithmetic index $\GA(\mathcal{G})$, Harmonic index $H(\mathcal{G})$ and Sum-Connectivity index  $\SCI(\mathcal{G})$ of $\mathcal{G}$  are defined as
\[
R(\mathcal{G})=\sum_{uv \in e(\mathcal{G})}\left(\deg_{\mathcal{G}}(u)\deg_{\mathcal{G}}(v)\right)^{\frac{-1}{2}}, \quad \ABC(\mathcal{G})=\sum_{uv\in e(\mathcal{G})}\left(\frac{\deg_{\mathcal{G}}(u)+\deg_{\mathcal{G}}(v)-2}{\deg_{\mathcal{G}}(u)\deg_{\mathcal{G}}(v)}\right)^{\frac{1}{2}}, 
\]
\[
\GA(\mathcal{G})=\sum_{uv\in e(\mathcal{G})} \frac{\sqrt{\deg_{\mathcal{G}}(u)\deg_{\mathcal{G}}(v)}}{\frac{1}{2}(\deg_{\mathcal{G}}(u)+\deg_{\mathcal{G}}(v))}, \quad H(\mathcal{G})=\sum_{uv\in e(\mathcal{G})}\frac{2}{\deg_{\mathcal{G}}(u)+\deg_{\mathcal{G}}(v)} 
\]
and
\[
\SCI(\mathcal{G})=\sum_{uv \in e(\mathcal{G})} \left(\deg_{\mathcal{G}}(u)+\deg_{\mathcal{G}}(v)\right)^{\frac{-1}{2}}.
\]
In this section, we compute the above mentioned topological indices of $\mathcal{B}(G)$ for the groups considered in Section 2. Since $\deg_{\mathcal{B}(G)}(x) = 1$ for any vertex $x \in G \times G$ we have the following result.  
\begin{lemma}\label{other_topological_index}
For any finite group $G$ we have 
  
\noindent  $R(\mathcal{B}(G))=\!\!\sum\limits_{H \in L(G)}\left(\deg_{\mathcal{B}(G)}(H)\right)^{\frac{1}{2}}$,
$\ABC(\mathcal{B}(G))=\!\!\sum\limits_{H \in L(G)}\left(\left(\deg_{\mathcal{B}(G)}(H)\right)^2-\deg_{\mathcal{B}(G)}(H)\right)^{\frac{1}{2}}$,
		 
$\GA(\mathcal{B}(G))=\sum\limits_{H \in L(G)} \frac{2\left(\deg_{\mathcal{B}(G)}(H)\right)^{\frac{3}{2}}}{(1+\deg_{\mathcal{B}(G)}(H))}$, \qquad
		 $H(\mathcal{B}(G))=\sum\limits_{H \in L(G)}\frac{2\deg_{\mathcal{B}(G)}(H)}{1+\deg_{\mathcal{B}(G)}(H)}$
		 
and	\quad	 $\SCI(\mathcal{B}(G))=\sum\limits_{H \in L(G)} \left(1+\deg_{\mathcal{B}(G)}(H)\right)^{\frac{-1}{2}}\deg_{\mathcal{B}(G)}(H)$.
\end{lemma}
While computing Zagreb indices of $\mathcal{B}(G)$ for various groups in Section 3, we have computed $\deg_{\mathcal{B}(G)}(H)$ for all $H \in L(G)$  (see the proofs of Theorem \ref{cyclic case} -- Theorem \ref{dycyclic-2}). Using those degrees of $H$ in $\mathcal{B}(G)$  and Lemma \ref{other_topological_index} we get the following results.
\begin{theorem}
Let $G$ be a cyclic group and $p$ be an odd prime. 
\begin{enumerate}
\item If $|G|=2p$ then 
 $R(\mathcal{B}(G))=(1+\sqrt{3})(1+\sqrt{p^2-1})$,
\begin{center}
$\ABC(\mathcal{B}(G))=\sqrt{6}+\sqrt{p^2-1}(\sqrt{p^2-2}+\sqrt{9p^2-12})$,
	 
	  $\GA(\mathcal{B}(G))=\frac{2+3\sqrt{3}}{2}+\sqrt{(p^2-1)^3}\left(\frac{2}{p^2}+\frac{6\sqrt{3}}{3p^2-2}\right)$,
\end{center}
\hspace{-1cm} $H(\mathcal{B}(G))=\frac{5}{2}+(p^2-1)\left(\frac{2}{p^2}+\frac{6}{3p^2-2}\right)$ and
		 $\SCI(\mathcal{B}(G))=\frac{1}{\sqrt{2}}+\frac{3}{2}+(p^2-1)\left(\frac{1}{p}+\frac{3}{\sqrt{3p^2-2}}\right)$.

\item If $|G|=2p^2$ then  $R(\mathcal{B}(G))=(1+\sqrt{3})(1+(p+1)\sqrt{p^2-1})$,

 $\ABC(\mathcal{B}(G))=\sqrt{6} \, +\sqrt{p^2-1} \, \left(\sqrt{p^2-2} \, +\sqrt{9p^2-12} \right.$

\qquad\qquad\qquad\qquad\qquad\qquad\qquad\qquad $+ \, \sqrt{p^6-p^4-p^2}+$  $\left. \sqrt{9p^6-9p^4-3p^2}\right)$,
\begin{center}
	 $\GA(\mathcal{B}(G))=\frac{2+3\sqrt{3}}{2}+\sqrt{(p^2-1)^3}\left(\frac{2}{p^2}+\frac{6\sqrt{3}}{3p^2-2}+\frac{2p^3}{p^4-p^2+1}+\frac{6p^3\sqrt{3}}{3p^4-3p^2+1}\right)$,

 $H(\mathcal{B}(G))=\frac{5}{2}+(p^2-1)\left(\frac{2}{p^2}+\frac{6}{3p^2-2}+\frac{2p^2}{p^4-p^2+1}+\frac{6p^2}{3p^4-3p^2+1}\right)$ and

 $\SCI(\mathcal{B}(G))=\frac{1}{\sqrt{2}}+\frac{3}{2}+(p^2-1)\left(\frac{1}{p}+\frac{3}{\sqrt{3p^2-2}}+\frac{p^2}{\sqrt{p^4-p^2+1}}+\frac{3p^2}{\sqrt{3p^4-3p^2+1}}\right)$.
\end{center}		
\item If $|G|=4p$ then 
 $R(\mathcal{B}(G))=(1+3\sqrt{3})(1+\sqrt{p^2-1})$,
\begin{center}
	 $\ABC(\mathcal{B}(G))=\sqrt{6}+\sqrt{132} +\sqrt{p^2-1}\left(\sqrt{p^2-2}+\sqrt{9p^2-12}+ \sqrt{144p^2-156}\right)$,

$\GA(\mathcal{B}(G))=\frac{26+135\sqrt{3}}{26}+\sqrt{(p^2-1)^3}\left(\frac{2}{p^2}+\frac{6\sqrt{3}}{3p^2-2}+\frac{48\sqrt{3}}{12p^2-11}\right)$,

$H(\mathcal{B}(G))=\frac{113}{26}+(p^2-1)\left(\frac{2}{p^2}+\frac{6}{3p^2-2}+\frac{24}{12p^2-11}\right)$ and

$\SCI(\mathcal{B}(G))=\frac{1}{\sqrt{2}}+\frac{3}{2}+\frac{12}{\sqrt{13}}+(p^2-1)\left(\frac{1}{p}+\frac{3}{\sqrt{3p^2-2}}+\frac{12}{\sqrt{12p^2-11}}\right)$.
\end{center}
\item If $|G|=4p^2$ then
$R(\mathcal{B}(G))=(1+3\sqrt{3})(1+(p+1)\sqrt{p^2-1})$,

\noindent  $\ABC(\mathcal{B}(G))=\sqrt{6}+\sqrt{132} +\sqrt{p^2-1}\left(\sqrt{p^2-2}+\sqrt{9p^2-12}+ \sqrt{144p^2-156}\right.$  

\qquad\quad\qquad$\left.+\sqrt{p^6-p^4-p^2}+\sqrt{9p^6-9p^4-3p^2}+\sqrt{144p^6-144p^4-12p^2}\right)$,

\noindent $\GA(\mathcal{B}(G))=\frac{26+135\sqrt{3}}{26}+\sqrt{(p^2-1)^3}\left(\frac{2}{p^2}+\frac{6\sqrt{3}}{3p^2-2}+\frac{48\sqrt{3}}{12p^2-11}\right.$  

\qquad\qquad\qquad\qquad\qquad\qquad\qquad\qquad$\left. +\frac{2p^3}{p^4-p^2+1}+ \frac{6p^3\sqrt{3}}{3p^4-3p^2+1}+\frac{48p^3\sqrt{3}}{12p^4-12p^2+1}\right)$,

\noindent$H(\mathcal{B}(G))=\frac{113}{26}+(p^2-1)\left(\frac{2}{p^2}+\frac{6}{3p^2-2}+\frac{24}{12p^2-11}+\frac{2p^2}{p^4-p^2+1}\right.$

\qquad\qquad\qquad\qquad\qquad\qquad\qquad\qquad\qquad\qquad$\left.+\frac{6p^2}{3p^4-3p^2+1} +\frac{24p^2}{12p^4-12p^2+1}\right)$ and

$\SCI(\mathcal{B}(G))=\frac{1}{\sqrt{2}}+\frac{3}{2}+\frac{12}{\sqrt{13}}+(p^2-1)\left(\frac{1}{p}+\frac{3}{\sqrt{3p^2-2}}+\frac{12}{\sqrt{12p^2-11}} \right.$ 

\qquad\qquad\qquad\qquad\qquad\qquad\qquad$\left. +\frac{p^2}{\sqrt{p^4-p^2+1}} +\frac{3p^2}{\sqrt{3p^4-3p^2+1}}+\frac{12p^2}{\sqrt{12p^4-12p^2+1}}\right)$.
\end{enumerate}    
\end{theorem}
\begin{theorem}
     If $G$ is a cyclic group of order $p^n$, where $p$ is a prime and $n \geq 1$, then 
\vspace*{.2cm} 
     \begin{center}
 $R(\mathcal{B}(G))=1+\sqrt{p^2-1}\left(\frac{p^n-1}{p-1}\right)$,
    \end{center}

$\ABC(\mathcal{B}(G))=\sqrt{p^2-1}\left(\sqrt{p^2-2}+p\sqrt{p^4-p^2-1}+p^2\sqrt{p^6-p^4-1}+\ldots\right.$

\qquad \qquad \qquad  \qquad \qquad \qquad \qquad \qquad \qquad \qquad \qquad \qquad $\left.+p^{n-1}\sqrt{p^{2n}-p^{2n-2}-1}\right)$,

\begin{center}
$\GA(\mathcal{B}(G))=1+\sqrt{(p^2-1)^3}\left(\frac{2}{p^2}+\frac{2p^3}{p^4-p^2+1}+\frac{2p^6}{p^6-p^4+1}+\ldots+\frac{2p^{3n-3}}{p^{2n}-p^{2n-2}+1}\right)$,

$H(\mathcal{B}(G))=1+(p^2-1)\left(\frac{2}{p^2}+\frac{2p^2}{p^4-p^2+1}+\frac{2p^4}{p^6-p^4+1}+\ldots+\frac{2p^{2n-2}}{p^{2n}-p^{2n-2}+1}\right)$ \qquad and

$\SCI(\mathcal{B}(G))=\frac{1}{\sqrt{2}}+(p^2-1)\left(\frac{1}{p}+\frac{p^2}{\sqrt{p^4-p^2+1}}+\frac{p^4}{\sqrt{p^6-p^4+1}}+\ldots+\frac{p^{2n-2}}{\sqrt{p^{2n}-p^{2n-2}+1}}\right)$.
\end{center}
\end{theorem}

\begin{theorem}
If $G=D_{2p}$, where $p$ is a prime, then
\begin{center}
 $R(\mathcal{B}(G))=1+p\sqrt{3}+\sqrt{p-1}(\sqrt{p+1}+\sqrt{3p})$,
 
$\ABC(\mathcal{B}(G))=p\sqrt{6}+\sqrt{p-1}\left(\sqrt{(p+1)(p^2-2)}+\sqrt{3p(3p^2-3p-1)}\right)$,

$\GA(\mathcal{B}(G))=\frac{2+3p\sqrt{3}}{2}+\sqrt{(p-1)^3}\left(\frac{2\sqrt{(p+1)^3}}{p^2}+\frac{6p\sqrt{3p}}{3p^2-3p+1}\right)$,

$H(\mathcal{B}(G))=\frac{2+3p}{2}+(p-1)\left(\frac{2p+2}{p^2}+\frac{6p}{3p^2-3p+1}\right)$ and

$\SCI(\mathcal{B}(G))=\frac{1}{\sqrt{2}}+\frac{3p}{2}+(p-1)\left(\frac{p+1}{p}+\frac{3p}{\sqrt{3p^2-3p+1}}\right)$.
\end{center}
\end{theorem}
\begin{theorem}
If $G=D_{2p^2}$, where $p$ is a prime, then
\begin{center}
$R(\mathcal{B}(G))=1+p^2\sqrt{3}+\sqrt{p-1}(\sqrt{p+1}(p+1)+2p\sqrt{3p})$,
\end{center}
 $\ABC(\mathcal{B}(G))=p^2\sqrt{6}+\sqrt{p^2-1}\left(\sqrt{p^2-2}+p\sqrt{p^4-p^2-1)}\right)$ 
 
\qquad \qquad\qquad\qquad\qquad\qquad $+p\sqrt{3p^2-3p}\left(\sqrt{3p^2-3p-1}+\sqrt{3p^4-3p^3-1}\right)$,

$\GA(\mathcal{B}(G))=\frac{2+3p^2\sqrt{3}}{2} + \sqrt{(p^2-1)^3}\left(\frac{2}{p^2} + \frac{2p^3}{p^4-p^2+1}\right)$ 

\qquad \qquad\qquad\qquad\qquad\qquad\qquad \qquad\qquad\qquad $+\sqrt{(p^2-p)^3}\left(\frac{2p}{3p^2-3p+1}+ \frac{2p^3}{3p^4-3p^3+1}\right)$,

\begin{center}
$H(\mathcal{B}(G))=\frac{2+3p^2}{2}+(p^2-1)\left(\frac{2}{p^2}+\frac{2p^2}{p^4-p^2+1}\right)+3p^2(p-1)\left(\frac{2}{3p^2-3p+1}+\frac{2p}{3p^4-3p^3+1}\right)$
\end{center}
and $\SCI(\mathcal{B}(G))=\frac{1}{\sqrt{2}}+\frac{3p^2}{2}+(p^2-1)\left(\frac{1}{p}+\frac{p^2}{\sqrt{p^4-p^2+1}}\right)$

\qquad \qquad\qquad\qquad\qquad\qquad\qquad \qquad\qquad$+3p^2(p-1)\left(\frac{1}{\sqrt{3p^2-3p+1}}+  \frac{p}{\sqrt{3p^4-3p^3+1}}\right)$.

\end{theorem}
\begin{theorem}
Consider the group $G=Q_{4p}$, where $p$ is a prime. 
\begin{enumerate}
\item If $p=2$ then 
 $R(\mathcal{B}(G))=1+7\sqrt{3}+2\sqrt{6}$, 
$\ABC(\mathcal{B}(G))= \sqrt{6}+3\sqrt{132}+\sqrt{552}$,
 $\GA(\mathcal{B}(G))=1+\frac{327\sqrt{3}}{26}+\frac{96\sqrt{6}}{25}$,
$H(\mathcal{B}(G))=\frac{4049}{650}$ and
$\SCI(\mathcal{B}(G))=\frac{1}{\sqrt{2}}+\frac{36}{\sqrt{13}}+\frac{63}{10}$.

\item If $p \geq 3$ then 
 $R(\mathcal{B}(G))=1+\sqrt{3}(1+2p)+\sqrt{p-1}(\sqrt{p+1}+\sqrt{3p+3}+\sqrt{12p})$,

\noindent $\ABC(\mathcal{B}(G))=\sqrt{6}+p\sqrt{132}$
 
\, $+\sqrt{p-1}\left(\sqrt{(p+1)(p^2-2)}+\sqrt{(p+1)(9p^2-12)} +\sqrt{12p(12p^2-12p-1)}\right)$,

\begin{center}
	 $\GA(\mathcal{B}(G))=\frac{2+3\sqrt{3}}{2}+\frac{24p\sqrt{12}}{13}+\sqrt{(p-1)^3}\left(\frac{2\sqrt{(p+1)^3}}{p^2}+\frac{2\sqrt{(3p+3)^3}}{3p^2-2}+\frac{2\sqrt{12p}}{12p^2-12p+1}\right)$,

$H(\mathcal{B}(G))=\frac{65+48p}{26}+(p-1)\left(\frac{2p+2}{p^2}+\frac{6p+6}{3p^2-2}+\frac{24p}{12p^2-12p+1}\right)$ and

 $\SCI(\mathcal{B}(G))=\frac{1}{\sqrt{2}}+\frac{3}{2}+\frac{12p}{\sqrt{13}}+(p-1)\left(\frac{p+1}{p}+\frac{3p+3}{\sqrt{3p^2-2}}+\frac{12p}{\sqrt{12p^2-12p+1}}\right)$.
\end{center}
\end{enumerate}
\end{theorem}
\begin{theorem}
Consider the group $G=Q_{4p^2}$, where $p$ is a prime. 
\begin{enumerate}
\item If $p=2$ then 
 $R(\mathcal{B}(G))=1+15\sqrt{3}+8\sqrt{6}$,

 $\ABC(\mathcal{B}(G))= \sqrt{6}+5\sqrt{132}+2\sqrt{552}+\sqrt{2256}+\sqrt{9120}$,
 
 $\GA(\mathcal{B}(G))=1+\frac{70830\sqrt{3}}{2548}+\frac{37824\sqrt{6}}{2425}$,
$H(\mathcal{B}(G))=\frac{60273113}{3089450}$ and

 $\SCI(\mathcal{B}(G))=\frac{1}{\sqrt{2}}+\frac{60}{\sqrt{13}}+\frac{96}{\sqrt{97}}+\frac{1257}{70}$.

\item If $p \geq 3$ then 

 $R(\mathcal{B}(G))=1+\sqrt{3}(1+2p^2)+\sqrt{p^2-1}(1+\sqrt{3} +p\sqrt{3})  $
 
\qquad\qquad\qquad\qquad $ +(p-1)\sqrt{12p^2-12p} + \sqrt{13p^4-12p^3+11p^2-12p})$,

\hspace{-1cm} $\ABC(\mathcal{B}(G))=\sqrt{6}+p^2\sqrt{132}+\sqrt{p^2-1}\left(\sqrt{p^2-2}+\sqrt{9p^2-12)} +p\sqrt{9p^4-9p^2-3}\right)$

\qquad \qquad $ +(p-1)\sqrt{(12p^2-12p)(12p^2-12p-1)}$
   
\qquad \qquad $+\sqrt{(13p^4-12p^3+11p^2-12p)(13p^4-12p^3+11p^2-12p-1)}$,

\noindent $\GA(\mathcal{B}(G))=\frac{2+3\sqrt{3}}{2}+\frac{24p^2\sqrt{12}}{13}+\sqrt{(p^2-1)^3}\left(\frac{2}{p^2}+\frac{6\sqrt{3}}{3p^2-2}+\frac{6p^3\sqrt{3}}{3p^4-3p^2+1}\right)$ 

\qquad \qquad\qquad \qquad\qquad \qquad\qquad$ +\frac{2(p-1)\sqrt{(12p^2-12p)^3}}{12p^2-12p+1}+\frac{2\sqrt{(13p^4-12p^3+11p^2-12p)^3}}{13p^4-12p^3+11p^2-12p+1}$,

\noindent $H(\mathcal{B}(G))=\frac{5}{2}+\frac{12p^2}{13}+(p^2-1)\left(\frac{2}{p^2}+\frac{6}{3p^2-2}+\frac{6p^2}{3p^4-3p^2+1}\right)$

\qquad \qquad\qquad \qquad\qquad \qquad\qquad \qquad$+\frac{2(p-1)(12p^2-12p)}{12p^2-12p+1}+\frac{2(13p^4-12p^3+11p^2-12p)}{13p^4-12p^3+11p^2-12p+1}$

\noindent and $\SCI(\mathcal{B}(G))=\frac{1}{\sqrt{2}}+\frac{3}{2}+\frac{12p^2}{\sqrt{13}}+(p^2-1)\left(\frac{1}{p}+\frac{3}{\sqrt{3p^2-2}}+\frac{3p^2}{\sqrt{3p^4-3p^2+1}}\right)$

\qquad \qquad\qquad \qquad\qquad \qquad\qquad \qquad$ +\frac{(p-1)(12p^2-12p)}{\sqrt{12p^2-12p+1}}+\frac{13p^4-12p^3+11p^2-12p}{\sqrt{13p^4-12p^3+11p^2-12p+1}}$.
\end{enumerate}
\end{theorem}

{\bf Acknowledgement.} The first author is thankful to Council of Scientific and Industrial Research  for the fellowship (File No. 09/0796(16521)/2023-EMR-I).

\section*{Declarations}
\begin{itemize}
	\item Funding: No funding was received by the authors.
	\item Conflict of interest: The authors declare that they have no conflict of interest.
	\item Availability of data and materials: No data was used in the preparation of this manuscript.
\end{itemize}

\end{document}